\documentclass[11pt,a4paper]{article}
\usepackage{amssymb,latexsym,amsmath}
\usepackage[english]{babel}
\newtheorem{theorem}{Theorem}
\newtheorem{lemma}{Lemma}
\newtheorem{definition}{Definition}
\newtheorem{notation}{Notation}

\newtheorem{example}{Example}
\begin{document}
\title{On the Number of $A$-mappings with Remainder}
\author{A.L. Yakymiv\\ Steklov Mathematical Institute of Russian Academy of Sciences\\
8 Gubkina St., Moscow 119991, Russia\\
arsen@mi.ras.ru}

\maketitle

\begin{abstract}
Mappings of a finite set into itself with restriction on the cycle lengths are considered (the so-called $A$-mappings). Asymptotics is given for the
number of these mappings with a power-law reduction of the remainder.

\end{abstract}

\section{Main Result}
Consider the sets $A$ having  positive densities $\varrho$ in the set of natural numbers, i.e., for which there exists
$$
\lim_{n\to\infty}\frac{|k: k\in A,\ k\leq n|}{n}=\varrho>0.\eqno(1.1)
$$
Put for $k\in N$
$$
p(k)=[s^k]\exp\left(\sum_{m\in
A}\frac{s^m}{m}\right),\ b(k)=kp(k),\ B(k)=\sum_{l=1}^k b(l).
\eqno(1.2)
$$
By $[s^k]g(s)$ we denote the coefficient for $s^k$ of the analytic function $g(s)$ in a neighborhood of zero.
Assume that for $k\to\infty$
$$
B(k)=ck^{\alpha}(1+O(k^{-\beta}))\eqno(1.3)
$$
for some positive constants $\alpha,\beta$ and $c$ (since this will not improve our estimates, we further assume that $\beta<1$). As it is shown in \cite[relation (20)]{Y13},  it follows from (1.1) that
$\alpha=\varrho+1$.

Let $X$ be a finite set containing $n$ elements. By $V_n(A)$ denote the set of mappings of $X$
into itself with contour sizes belonging to the set $A$. Such mappings were introduced by V.N. Sachkov 
\cite{Sach72} (1972). By $\lambda_n(A)$, denote the number of cyclic elements of a random mapping
having uniform distribution on the set $V _{n}(A)$.

Put for $\mu>-1$
$$
I_\mu=\int_0^\infty
x^\mu\exp\left(-\frac{x^2}{2}\right)\,dx.\eqno(1.4)
$$
\begin{theorem}
Let (1.3) holds. Then
$$
|V_{n}(A)|=c(1+\varrho)n^{n-(1-\varrho)/2}(I_\varrho+O(n^{-\beta/2})),\eqno(1.5)
$$
as $n\to\infty$. In addition, for an arbitrary fixed $z>0$
$$
\mathsf{P}\left\{\frac{\lambda_n}{\sqrt{n}}\leq z\right\}=
\frac{1}{I_\varrho}\int_0^zx^\varrho\exp\left(-\frac{x^2}{2}\right)\,dx
+O(n^{-\beta/2}).\eqno(1.6)
$$
\end{theorem}
\begin{notation}
The first term in (1.5) is obtained in \cite{Y13} under assumption (1.1).
\end{notation}
Let us give some examples of sets $A$ for which relation (1.3) holds.
\begin{definition}
We say that the set $A$ belongs to
class $F_{1}$ iff
$A=\bigcup_{i=1}^{M}A_{i}$, where $M\in N,\ A_{i}=\{m\in N:\
m=a_{i}k+b_{i},\ k=0,1,2,\dots \}$ and
integers $a_{i}>1,\ 1\leq b_{i}\leq a_{i}-1,\ (a_{i},b_{i})=1$, and
the progressions $A_{i}$ and $A_{j}$ do not intersect for $i\neq j$.
\end{definition}

\begin{definition}
We say that the set $A$ belongs to
class $F_{2}$ if $A=\{m\in N:\
m/k_{i}\notin N,\ i=1,\dots ,s\}$ for some $s\in N$ and $k_{1},\dots
,k_{s}\in N$ such that
$k_{i}\geq 2,\ i=1,\dots ,s$ and $(k_{i},k_{j})=1$ for $i\neq j$.
\end{definition}
\begin{example}
Let the set $A\in F_{1}\cup F_{2}$. In \cite{P97}, it is proved that for some positive constants $C$ and $\beta$ (depending on $A$) for
$n\rightarrow\infty$
$$
p(n)=Cn^{\varrho-1}\left(1+O\left(n^{-\beta}\right)\right).\eqno(1.7)
$$
\end{example}
It follows from (1.7) that for $A\in F_{1}\cup F_{2}$ relation (1.3) holds with $\alpha=\varrho+1$ and $c=C/(\varrho+1)$.
In fact, we obtain from (1.7) that
$$
B(n)=\sum_1^nkp(k)=C\sum_1^nk^{\varrho}+O(1)\sum_1^nk^{\varrho-\beta}
$$
$$
=\frac{Cn^{\varrho+1}}{\varrho+1}+O(n^\varrho)+O(1)n^{\varrho-\beta+1}=cn^{\alpha}(1+O(n^{-\beta})),
$$
since
$$
\sum_1^nk^{\varrho}=\frac{n^{\varrho+1}}{\varrho+1}+O(n^\varrho)\eqno(1.8)
$$
and $\beta<1$. The following example shows that (1.3) does not imply (1.7).
\begin{example}
Let $A=aN$, where the natural number $a>1$ is fixed. Then the set $A$ satisfies (1.3).
\end{example}
However, in this case, relation (1.7) does not hold since $p(n)=0$ for $n\in N\backslash A$. The assertion of this example is proven in Section 2.

\section{Auxiliary propositions}
Set for $k,n\in N$ and $k\leq n$
$$
a(k,n)=\left(1-\frac{1}{n}\right)
\cdots\left(1-\frac{k-1}{n}\right). \eqno(2.1)
$$
\begin{lemma}
For all $k\in(1,n]$, the following inequality holds:
$$
a(k,n)\leq\sqrt{e}\exp\left(-\frac{k^2}{2n}\right).
$$
\end{lemma}
Lemma 2 follows from inequality (14) \cite{Y9}.
\begin{lemma}
For $n\to\infty$
$$
a(k,n)=\exp\left(-\frac{k^2}{2n}\right)\left(1+O\left(\frac{k^3}{n^2}\right)\right).
$$
uniformly by $k\in[1,n^{2/3}]\cap N$.
\end{lemma}
This lemma follows from inequality (13) \cite{Y9}. Set for $\mu>-1$
$$
\Psi_{\mu}(n)=\sum_{r(n)+1}^{s(n)}\exp\left(-\frac{k^2}{2n}\right)k^{\mu},\eqno(2.2)
$$
$$
\Sigma_\mu(n)=\sum_{r(n)+1}^{s(n)}\frac{1}{\sqrt{n}}\left(\frac{k}{\sqrt{n}}\right)^{\mu}\exp\left(-\frac{k^2}{2n}\right)
=n^{-(1+\mu)/2}\Psi_{\mu}(n)\eqno(2.3)
$$
and
$$
I_\mu(n)=\int_{(r(n)+1)/\sqrt{n}}^{(s(n)+1)/\sqrt{n}}x^\mu\exp\left(-\frac{x^2}{2}\right)\,dx,\eqno(2.4)
$$
where the natural numbers $r(n)<s(n)$ will be chosen later.
\begin{lemma}
$$
-\sum_{r(n)+1}^{s(n)}\frac{1}{\sqrt{n}}\left(\frac{k}{\sqrt{n}}\right)^{\mu}\exp\left(-\frac{k^2}{2n}\right)\left(\left(1+\frac{1}{k}\right)^\mu-1\right)
\leq \Sigma_\mu(n)-I_\mu(n)\leq\frac{2}{\sqrt{n}}\Sigma_\mu(n)
$$
\end{lemma}
This lemma follows from the chain of inequalities (20) and (21) of the article \cite{Y9}.
\begin{lemma}\cite[lemma 3]{P97}
For an arbitrary fixed $\lambda>0$ we have:
$$
[s^n](1-s)^{-\lambda}=\frac{1}{\Gamma(\lambda)}n^{\lambda-1}\left(1+O\left(\frac{1}{n}\right)\right).
$$
\end{lemma}
Proof of the assertion of example 2. We have:
$$
\sum_{m\in A}\frac{s^m}{m}=\sum_{k\in N}\frac{s^{ak}}{ak}=\frac{1}{a}\sum_{k\in N}\frac{s^{ak}}{k}
=-\frac{1}{a}\ln(1-s^a),
$$
since
$$
\sum_{k\in N}\frac{x^{k}}{k}=-\ln(1-x).
$$
for $x\in(0,1)$. Thus
$$
\exp\left(\sum_{m\in A}\frac{s^m}{m}\right)=(1-s^a)^{-1/a},
$$
whence
$$
p(an)=[s^{an}]\left(\frac{1}{1-s^a}\right)^{1/a}=[s^{n}]\left(\frac{1}{1-s}\right)^{1/a}
=\frac{n^{1/a-1}}{\Gamma(1/a)}\left(1+O\left(\frac{1}{n}\right)\right)
$$
according to Lemma 4.  We obtain from here that for $m=[k/a]$
$$
B(k)=\sum_{l=1}^k lp(l)=\sum_{n=1}^{m}anp(an)
=\frac{a}{\Gamma(1/a)}\sum_{n=1}^{m}n^{1/a}\left(1+O\left(\frac{1}{n}\right)\right).\eqno(2.5)
$$
By virtue of (1.8) we have:
$$
\sum_{n=1}^{m}n^{1/a}=\frac{m^{1/a+1}}{1/a+1}+O(m^{1/a})=\frac{m^{1/a+1}}{1/a+1}\left(1+O\left(\frac{1}{m}\right)\right)
$$
$$
=\frac{m^{1/a+1}}{1/a+1}\left(1+O\left(\frac{1}{k}\right)\right).\eqno(2.6)
$$
From (2.5) and (2.6) we obtain that
$$
B(k)=\frac{am^{1/a+1}}{(1/a+1)\Gamma(1/a)}\left(1+O\left(\frac{1}{k}\right)\right)
$$
$$
=k^{1/a+1}\frac{a^{-1/a}}{(1/a+1)\Gamma(1/a)}\left(1+O\left(\frac{1}{k}\right)\right),\eqno(2.7)
$$
since
$$
m^{1/a+1}=[k/a]^{1/a+1}=(k/a+O(1))^{1/a+1}
$$
$$
=\left(k/a\left(1+O\left(\frac{1}{k}\right)\right)\right)^{1/a+1}=(k/a)^{1/a+1}\left(1+O\left(\frac{1}{k}\right)\right).
$$
It follows from (2.7) that (1.3) holds for the set $A=aN$  with $\alpha=1/a+1,\ \beta=1$ and $c=(a^{1/a}(1/a+1)\Gamma(1/a))^{-1}$.

\section{Proof of Theorem 1}
As it is shown in the paper \cite[relations (10)-(14)]{Y13},
$$
|V_n(A)|
=n^{n-1}\sum_{k=1}^na(k,n)b(k)
$$
$$
=n^{n-1}\left(\frac{(n-1)!}{n^n}B(n)+\frac{1}{n}\sum_{k=1}^{n-1}a(k,n)kB(k)\right).\eqno(3.1)
$$
It follows from (1.3) and the equality $\alpha=\varrho+1$ that
$$
\frac{(n-1)!}{n^n}B(n)=O(n^\varrho)\frac{n!}{n^n}=O(1)e^{-n}n^{\varrho+1/2}
$$
Thus
$$
|V_n(A)|=n^{n-2}\left(\sum_{k=1}^{n-1}a(k,n)kB(k)+O\left(e^{-n}n^{\varrho+3/2}\right)\right).\eqno(3.2)
$$
Further,
$$
\sum_{k=1}^{n-1}a(k,n)kB(k)=\left(\sum_{1}^{r(n)}+\sum_{r(n)+1}^{s(n)}+\sum_{s(n)+1}^{n-1}\right)a(k,n)kB(k)
$$
$$
=J_1+J_2+J_3.\eqno(3.3)
$$
According to (2.1) $a(k,n)<1$. Therefore, taking into account (1.3), we obtain:
$$
J_1=\sum_{1}^{r(n)}a(k,n)kB(k)<\sum_{1}^{r(n)}kB(k)=O(1)\sum_{1}^{r(n)}k^{\varrho+2}=O(1)r(n)^{\varrho+3}\eqno(3.4)
$$
It follows from (1.3) that
$$
J_2=\sum_{r(n)+1}^{s(n)}a(k,n)kB(k)=c\sum_{r(n)+1}^{s(n)}a(k,n)k^{2+\varrho}+O(1)\sum_{r(n)+1}^{s(n)}a(k,n)k^{2+\varrho-\beta}.
$$
Hence, we have by Lemma 1 that
$$
J_2=c\sum_{r(n)+1}^{s(n)}\exp\left(-\frac{k^2}{2n}\right)k^{2+\varrho}+O(1)\sum_{r(n)+1}^{s(n)}\exp\left(-\frac{k^2}{2n}\right)k^{2+\varrho}\left(k^{-\beta}+ \frac{k^3}{n^2}\right).\eqno(3.5)
$$
Note that the inequality
$$
\frac{k^3}{n^2}\leq k^{-\beta}
$$
holds  iff
$$
k^{3+\beta}\leq n^2
$$
or
$$
k\leq\left[n^{\frac{2}{3+\beta}}\right]\stackrel{def}=s(n).
$$
Therefore, setting $r(n)=[\ln n]$,  we obtain from (3.5) that
$$
J_2=c\sum_{r(n)+1}^{s(n)}\exp\left(-\frac{k^2}{2n}\right)k^{2+\varrho}+O(1)\sum_{r(n)+1}^{s(n)}\exp\left(-\frac{k^2}{2n}\right)k^{2+\varrho-\beta}
$$
$$
=c\Psi_{2+\varrho}(n)+O(1)\Psi_{2+\varrho-\beta}(n)
$$
$$
=cn^{(3+\varrho)/2}\Sigma_{2+\varrho}(n)+O(1)n^{(3+\varrho-\beta)/2}\Sigma_{2+\varrho-\beta}(n)
$$
$$
=n^{(3+\varrho)/2}(c\Sigma_{2+\varrho}(n)+O(n^{-\beta/2})).\eqno(3.6)
$$
According to (1.3),
$$
J_3=\sum_{s(n)+1}^{n}a(k,n)kB(k)=O(1)\sum_{s(n)+1}^{n}a(k,n)k^{\varrho+2}.
$$
Thus, we have by Lemma 1 that
$$
J_3=O(1)\sum_{s(n)+1}^{n}k^{\varrho+2}\exp\left(-\frac{k^2}{2n}\right)
$$
$$
=O\left(n^{(3+\varrho)/2}\right)\sum_{s(n)+1}^{n}\frac{1}{\sqrt{n}}\left(\frac{k}{\sqrt{n}}\right)^{\varrho+2}\exp\left(-\frac{k^2}{2n}\right)
$$
$$
=O\left(n^{(3+\varrho)/2}\right)\int_{s(n)/\sqrt{n}}^\infty x^{\varrho+2}\exp\left(-\frac{x^2}{2}\right)\,dx.\eqno(3.7)
$$
Note that
$$
\nu\stackrel{def}=\frac{2}{3+\beta}-\frac{1}{2}=\frac{1-\beta}{2(3+\beta)}>0,
$$
since $\beta<1$. Hence, we obtain from (3.5) that
$$
J_3=O\left(n^{(3+\varrho)/2}\right)\int_{n^\nu}^\infty x^{\varrho+2}\exp\left(-\frac{x^2}{2}\right)\,dx
$$
$$
=O\left(n^{(3+\varrho)/2}\right)O(n^{-\beta/2})=O\left(n^{(3+\varrho-\beta)/2}\right).\eqno(3.8)
$$
Note that  
$$
\sum_{r(n)+1}^{s(n)}\frac{1}{\sqrt{n}}\left(\frac{k}{\sqrt{n}}\right)^{\mu}\exp\left(-\frac{k^2}{2n}\right)\left(\left(1+\frac{1}{k}\right)^\mu-1\right)
$$
$$
\sim\varrho\sum_{r(n)+1}^{s(n)}\frac{1}{\sqrt{n}}\left(\frac{k}{\sqrt{n}}\right)^{\varrho}\exp\left(-\frac{k^2}{2n}\right)\frac{1}{k}
=\frac{\varrho}{\sqrt{n}}\Sigma_{\mu-1}(n)=O\left(\frac{1}{\sqrt{n}}\right).
$$
Therefore, we obtain from Lemma 3 that
$$
\Sigma_{2+\varrho}(n)=I_{2+\varrho}(n)+O\left(\frac{1}{\sqrt{n}}\right).\eqno(3.9)
$$
It follows from (3.7) and (3.9) that
$$
J_2=n^{(3+\varrho)/2}\left(cI_{2+\varrho}(n)+O(1)n^{-\beta/2}+O\left(n^{-1/2}\right)\right)
$$
$$
n^{(3+\varrho)/2}(cI_{2+\varrho}(n)+O(n^{-\beta/2}))\eqno(3.10)
$$
since $\beta<1$. We deduce from (3.3), (3.4), (3.8), (3.9) and (3.10) that
$$
J=J_1+J_2+J_3
$$
$$
=O(\ln^{\varrho+3}n)+n^{(3+\varrho)/2}(cI_{2+\varrho}(n)+O(n^{-\beta/2}))+O\left(n^{(3+\varrho-\beta)/2}\right)
$$
$$
=n^{(3+\varrho)/2}(cI_{2+\varrho}(n)+O(n^{-\beta/2})).
$$
Thus, it follows from (3.2) and (3.3) that
$$
|V_n(A)|=n^{n-2}\left(J+O\left(e^{-n}n^{\varrho+3/2}\right)\right)
$$
$$
=n^{n-2}\left(n^{(3+\varrho)/2}(cI_{2+\varrho}(n)+O(n^{-\beta/2}))+O\left(e^{-n}n^{\varrho+3/2}\right)\right)=
$$
$$
=n^{(3+\varrho)/2}n^{n-2}\left(cI_{2+\varrho}(n)+O(n^{-\beta/2}))\right)
$$
$$
=n^{n-(1-\varrho)/2}\left(cI_{2+\varrho}(n)+O(n^{-\beta/2}))\right).\eqno(3.11)
$$
Note that
$$
I_{2+\varrho}(n)=I_{2+\varrho}-\left(\int_0^{(\ln n+1)/\sqrt{n}}+\int_{[n^\nu]}^\infty\right)x^{2+\varrho}\exp\left(-\frac{x^2}{2}\right)\,dx
$$
$$
=I_{2+\varrho}+O(n^{-\beta})=(\varrho+1)I_{\varrho}+O(n^{-\beta}).\eqno(3.12)
$$
Relation (1.5) follows from (3.11) and (3.12).
Further, according to formula (26) from paper \cite{Y13}, for
for an arbitrary fixed positive $z$ we have:
$$
\mathsf{P}\left\{\frac{\lambda_n}{\sqrt{n}}\leq z\right\}=
\frac{\sum_{k=1}^{\left[z\sqrt{n}\right]}a(k,n)b(k)}{\sum_{k=1}^n
a(k,n)b(k)}.\eqno(3.13)
$$
We obtain from (3.1) and (1.5) that
$$
\sum_{k=1}^n
a(k,n)b(k)=n^{1-n}|V_n(A)|\sim n^{(1+\varrho)/2}c(1+\varrho)\left(I_\varrho
+O(n^{-\beta/2})\right).\eqno(3.14)
$$
To find the asymptotic formula for the numerator of the fraction from (3.13), we carry out the same reasoning as in the proof of (1.5). As a result, we get that for $m=\left[z\sqrt{n}\right]$
$$
\sum_{k=1}^{m}a(k,n)b(k)=a(m,n)B(m)-\sum_{k=1}^{m-1}(a(k+1,n)-a(k,n))B(k)=
$$
$$
=a(m,n)B(m)+\frac{1}{n}\sum_{k=1}^{m-1}a(k,n)kB(k)
$$
$$
\sim n^{(1+\varrho)/2}c(1+\varrho)\left(\int_0^zx^\varrho\exp\left(-\frac{x^2}{2}\right)\,dx
+O(n^{-\beta/2})\right).\eqno(3.15)
$$
Relation (1.6) follows from (3.13), (3.14) and (3.15). The proof of the theorem is complete.

\end{document}